\documentclass[12pt]{article}
\usepackage{a4,amsmath,amsthm,amssymb,amsfonts}
\usepackage{verbatim}
\theoremstyle{plain}
\newtheorem{thm}{Theorem}
\newtheorem{cor}{Corollary}
\newtheorem{ack}{Acknowledgement \hskip -6 pt}

\newtheorem{lem}{Lemma}
\newtheorem{prop}{Proposition}

\theoremstyle{definition}   
\newtheorem{rem}{Remark}

\newtheorem{defn}{Definition}
\newtheorem{example}{Example}

\renewcommand{\Re}{\operatorname{Re}}
\title{On Factorisations of Matrices and Abelian Groups.}
\author{Johan Andersson\footnote{Department of Mathematics, Stockholm University, SE-106 91 Stockholm, Sweden. Email : johana@math.su.se} \, and Gautami Bhowmik\footnote{Universit\'e de Lille
    1, UMR CNRS 8524, 59655 Villeneuve d'Ascq Cedex, France. Email : bhowmik@math.univ-lille1.fr}}

\date{\today}

\begin{document}

\maketitle

\begin{abstract} 
  We establish correspondances between factorisations of finite abelian groups
  ( direct factors, 
  unitary factors, non isomorphic 
     subgroup classes ) and factorisations of integer matrices.
    We then study counting functions associated to these factorisations and find average orders.
\end{abstract}


\section{ Introduction}
 In this paper we study the correspondance between different factorisations of
finite abelian groups and that of integer matrices. In particular, we are
interested in the matter of enumeration.  This is in line with the
questions on direct and unitary factors of abelian groups first
raised by Cohen \cite{C} in 1960 and subgroups of abelian groups studied since
long
in different contexts
(see, eg, Butler \cite{Butler}, Goldfeld-Lubotzky-Pyber \cite{GLP} or Bhowmik-Ramare \cite{BR}).

We consider the decomposition of a 
finite abelian group $G$ 
as a direct sum
$$
G\cong 
\mathbb Z/n_1\mathbb Z\oplus\cdots\oplus
\mathbb Z/n_1n_2\mathbb Z
\oplus\cdots\oplus\mathbb Z/n_1n_2\cdots n_r\mathbb Z
$$
where  $n_i$ are divisors of the {\it order} $n=n_1^r n_2^{r-1}\cdots n_r$ 
and $r$ the {\it rank} of the group.
In terms of integer matrices this group
can be
 associated to  a canonical
representative of a double coset $$M(G)=diag[n_1,n_1n_2,\cdots,n_1n_2\cdots n_r].$$ 
The matrix $M(G)$ is said to be
in {\it Smith Normal Form} (SNF) or sometimes an elementary divisor and its determinant equals the order
of the group $G$. 
 It is well 
known that the number of non-isomorphic abelian groups of order upto $x$
 is asymptotic to $Lx$ \cite{ES}  where $L$ is a constant.  If the rank is bounded
by $r$, this corresponds to the number of
$r\times r$ SNF matrices of determinant upto $x$ and then the constant depends on $r$ \cite{B}.

In 1960 E. Cohen \cite{C} defined the factorization of finite abelian groups
as a decomposition in formal direct factors. He introduced the
zeta function associated to the direct factorization, $G=H\oplus K$,
 and obtained
that the number of direct factors, on an average, is $\log x$. In Section 1
we interpret direct factorisations of $G$ in terms of decompositions into SNF
factors of the `conjugate' matrix of $G$. The idea of conjugate matrices
comes from considering  the type of a group, in the sense of Macdonald \cite{Md}. Thus we can treat the truncated zeta function
$$\zeta^2(s)\zeta^2(2s)\zeta^2(3s)\cdots\zeta^2(rs)$$ as a counting function 
of certain direct factors of abelian groups.  The case of unitary factorisations
is treated similarly.

Later, \cite{BR} it was proved 
that a subgroup of $G$ corresponds
to a divisor of $M(G)$ which is a canonical representative
of a right coset, for example a matrix in {\it Hermite Normal Form} (HNF).  
Using the zeta function associated to subgroups
it was found that on an average there were many more subgroups than
direct factors. The number of subgroups of an abelian group
is essentially $x^{([r^2/4]+1)/r -1}$, which is heavily dependent on the rank. 

In Section 2 we consider classes of subgroups under isomorphism. In terms of
matrices we are interested in divisors of $M(G)$ which are in SNF. We now  investigate whether
on an average the number of non-isomorphic subgroups  depends on the rank of the group.

To do so we introduce the zeta function
 $$Z_r(s)=\sum_{n=0}^\infty A_r(n) n^{-s}$$
 where $A_r(n)$ denotes the number of subgroup classes (under isomorphism) of
 abelian groups of order  $n$ rank atmost $r$ 
 and  use this function to get the mean value
  $$
  \sum_{n=1}^x A_r(n) \sim C_r x  \log x
  $$
for positive integers $r$. This   shows that on an average the number of non-isomorphic subgroups is
about the same as the number of direct factors whereas the number of all
subgroups is much larger as the rank increases. 
The same average order holds true when we consider subgroups of arbitrary rank.

 While $Z_r(s)$ in
general have no global functional equation we will prove 
that for all positive integers $r$ it has  a local functional equation 
\begin{gather*}
  Z_r(s)=\prod_{p \text{ prime}} F_r(p^{-s}), \\ \intertext{with}
 F_r(x)=x^{-r(r+3)/2} F_r \left(\frac 1 x \right).
\end{gather*}
This follows from a recursion formula that we develop for the functions $F_r(x)$.
 We  include some values of the generating function by using this
recursion formula, the first few being
 $$Z_1(s)=\zeta^2(s), \quad  Z_2(s)=\frac{\zeta^2(s)
\zeta^3(2s)}{\zeta(3s)}\quad\text{and} \quad Z_3(s)= \frac{\zeta^2(s) \zeta^3(2s) \zeta^3(3s)}{ \zeta(4s)} G(s), $$
where $G(s)$ is a Dirichlet series absolutely convergent for $\Re(s)>1/5$. 

For the limit case $Z(s)=\lim_{r \to\infty} Z_r(s)$ we show that while $Z(s)$ does have a meromorphic continuation to $\Re(s)>0$,
 it does not have a meromorphic continuation beyond the imaginary axis. We 
tend to believe that the same is true for $Z_r(s)$ when $r \geq 3$.

\begin{ack}   The second author thanks Radoslav
Dimitri\c c for useful discussions.
\end{ack}

\section{ Direct Factors.}

In this section we discuss  the relationship between direct factors
of abelian groups and factorisation of matrices.
We recall a definition.

\begin{defn} An abelian group $G$ has a (formal) {\it  direct factorisation} into
$H$ and $K$, noted as $G=H\oplus K$ if $G= H + K$, $H$ and $K$ being subgroups
of $G$ whose only
intersection as subgroups is trivial. \end{defn}

Formal direct factors of abelian groups were studied in the
context of arithmetic functions of abelian groups and
Cohen \cite{C} defined direct Dirichlet convolution of abelian groups
in terms of this decomposition, i.e.
$$   
(\chi_1*\chi_2)(G)=\sum_{G=H\oplus K}\chi_1(H)\chi_2(K),
$$
 where $\chi_i$ are arithmetical functions of abelian groups, in an attempt to
 generalise the factorisation of integers.

With no loss of generalisation, we can consider $p$-groups, i.e.
$$
G_p\cong 
\mathbb Z/p^{a_1}\mathbb Z
\oplus\mathbb Z/p^{a_2}\mathbb Z\oplus\cdots
\oplus
\mathbb Z/p^{a_r}\mathbb Z,
$$
where the order of $G_p$ is $p^{a_1+a_2+\cdots+ a_r}$, for a prime number $p$
and $0\le a_1\le a_2\le\cdots\le a_r$.

 The invariants of $G$ is given by the
partition $\lambda=(a_1,a_2,\cdots, a_r)$, which, in combinatorics, is called 
the {\sl type} of $G$. 
We call the group whose type $\lambda'$
is the conjugate of $\lambda$, the {\sl conjugate of G}
and denote it by $G'$. 
We now express the direct factorisation above as a matrix factorisation.

\begin{defn}The {\sl SNF factorisation} of a SNF
matrix $S$ is a decomposition into SNF matrices whose product
equals $S$. \end{defn}

We now establish the following bijection
\begin{prop} \quad If $M(G')=M(H')M(K')$ is a SNF factorisation, then
$G=H\oplus K$ is a direct factorisation and conversely. \end{prop}

\begin{proof}
Let $G$ be a group of type 
\begin{align*}
\lambda &=
(\underbrace {f_1}_{g_1\ {\text {times}}},\quad
\underbrace{f_1+f_2}_{g_2\ {\text {times}}},\dots,
\underbrace{f_1+f_2+\dots+f_{m}}_{g_{m}\ {\text {times}}}).\\ \intertext{Then} \lambda'&=
(\underbrace {g_m}_{f_m\ {\text {times}} },\quad
\underbrace{g_m+g_{m-1}}_{f_{m-1}\ {\text {times}}},\dots,
\underbrace{g_m+g_{m-1}+\dots+g_{1}}_{f_{1}\ {\text {times}}}). \end{align*}
Writing $g_i=h_i+k_i$, sum of non-negative integers, we get types
of groups $H'$ and $K'$ such that
$M(G')=M(H')M(K')$ is a SNF factorisation, i.e.
$$
M(H')=diag[\underbrace{p^{h_m},\dots,p^{h_m}}_{f_m\ {\text {times}}}
,\underbrace{p^{h_m+h_{m-1}},\dots}_{f_{m-1}\ {\text {times}}},
\dots,\underbrace{p^{h_m+h_{m-1}+\dots+h_1}}
_{f_1\ {\text {times}}}]  \label{500}
$$
etc. whereas 
$$
H\cong \oplus_{h_1}\mathbb Z/p^{f_1}\mathbb Z
\oplus_{h_2}\mathbb Z/p^{f_1+f_2}\mathbb Z
\dots\oplus_{h_m}\mathbb Z/p^{f_1+f_2+\dots+f_m}\mathbb Z
$$
and
$$
K\cong \oplus_{k_1}\mathbb Z/p^{f_1}\mathbb Z
\oplus_{k_2}\mathbb Z/p^{f_1+f_2}\mathbb Z
\dots\oplus_{k_m}\mathbb Z/p^{f_1+f_2+\dots+f_m}\mathbb Z,
$$
where $\oplus_{t}$ before a summand means that it is repeated $t$ times. 

Thus the type of $H$ is
\begin{gather*}(\underbrace {f_1}_{h_1\ {\text {times}}},\quad
\underbrace{f_1+f_2}_{h_2\ {\text {times}}},\dots,
\underbrace{f_1+f_2+\dots+f_{m}}_{h_{m}\ {\text {times}}}) \\ \intertext{
and that of $K$ is}
(\underbrace {f_1}_{k_1\ {\text {times}}},\quad
\underbrace{f_1+f_2}_{k_2\ {\text {times}}},\dots,
\underbrace{f_1+f_2+\dots+f_{m}}_{k_{m}\ {\text {times}}}), \end{gather*}
 and the relation $G\cong H \oplus K$ is satisfied. 

The converse is proved by taking $G\cong H \oplus K$ as above and retracing
the steps.
\end{proof}

\begin{example}  Let $G\cong \mathbb Z/p^{2}\mathbb Z
\oplus\mathbb Z/p^{2}\mathbb Z
\oplus
\mathbb Z/p^{4}\mathbb Z$. Here $\lambda=(2,2,4)$, thus the conjugate 
$\lambda'=(1,1,3,3)$. The following non-trivial SNF factorisations of $M(G')$, i.e.
$$ M(G')=
\begin{pmatrix}
1      &0          &0     &0                 \cr
0      &1          &0     &0                 \cr
0      &0          &p     &0                 \cr
0      &0          &0     &p                 \cr
\end{pmatrix}
\begin{pmatrix}
p      &0      &0     &0       \cr
0      &p      & 0    &0       \cr
0      &0      &p^2   &0       \cr
0      &0      &0     &p^{2}   \cr
\end{pmatrix}
=\begin{pmatrix}
1      &0          &0     &0                 \cr
0      &1          &0     &0                 \cr
0      &0          &p^2   &0                 \cr
0      &0          &0     &p^2               \cr
\end{pmatrix}
\begin{pmatrix}
p      &0      &0     &0     \cr
0      &p      & 0    &0     \cr
0      &0      &p     &0     \cr
0      &0      &0     &p     \cr
\end{pmatrix}
$$
 give 
$H\cong \mathbb Z/p^{2}\mathbb Z, 
K\cong \mathbb Z/p^{2}\mathbb Z
\oplus\mathbb Z/p^{4}\mathbb Z$ and $H\cong \mathbb Z/p^{2}\mathbb Z
\oplus\mathbb Z/p^{2}\mathbb Z, K\cong \mathbb Z/p^{4}\mathbb Z.$ Obviously, changing the
order of the matrices would interchange $H$ and $K$.
\end{example}
\medskip

The generating function for the number of direct factors of a finite
abelian group is given by 
$$ 
{\mathcal D}(s)=\zeta^2(s)\zeta^2(2s)\zeta^2(3s)\cdots
$$
as was first done by Cohen \cite{C} and subsequently studied by many, for
example Kr\" atzel \cite{K}, Menzer \cite{M}, Wu \cite{W} or Knopfmacher \cite{Knopfmacher}.

A truncated form of this function, 
$$ 
{\mathcal D}_r(s)=\zeta^2(s)\zeta^2(2s)\zeta^2(3s)\cdots\zeta^2(rs)
$$
occurs as the Dirichlet series associated to the SNF factorisation
of an $r\ \times r $ integer matrix, see, for example \cite{B}, which we can
now interpret in the context of abelian groups: 
\begin{thm}
The function
$\zeta^2(s)\zeta^2(2s)\zeta^2(3s)\cdots\zeta^2(rs)$ generates the direct
factors of finite abelian groups
whose conjugates have rank atmost $r$. 
\end{thm}
This implies the following result on average orders:
\begin{cor} 
 The number of  direct
factors of non-isomorphic  abelian groups of order atmost $x$ whose conjugates
have rank atmost $r$ is asymptotic to $A_r x\log x$, where $A_r$ is a constant that depends on
$r$.
\end{cor}

\bigskip
\section{ Unitary Factors.} 

Unitary factors of finite abelian groups were also studied by
Cohen as a generalisation of the corresponding idea on integers, 
where one encounters the function 
$$t(n)=\sum_{d\mid n,\ gcd(d,n/d)=1}1.$$

\begin{defn} An abelian group $G$ has  unitary factors
$H$ and $K$ if $G\cong H \oplus K$ is a direct factorisation $H$ and $K$ have
no further common direct factor. \end{defn}

As above, we see that in terms of matrices this corresponds to block SNF factorisations, i.e. 
for $M(G')$ as above, $M(H')$ is a unitary factor if $h_i$ is either 
$0$ or $g_i$ for each $i$. Thus the number of unitary factors of $G$
is $2^m$.
\medskip 
\begin{example} With $G$ as in Example 1, the only non-trivial unitary 
factorisation $H\cong \mathbb Z/p^{2}\mathbb Z
\oplus\mathbb Z/p^{2}\mathbb Z, K\cong \mathbb Z/p^{4}\mathbb Z$ 
corresponds to the SNF factorisation 
$$ M(G')= \begin{pmatrix}
pI_2      & 0                         \cr
0         &p^3I_2                    \cr
\end{pmatrix}
=\begin{pmatrix}
I_2      & 0                         \cr
0         &p^2I_2                    \cr
\end{pmatrix}
\begin{pmatrix}
pI_2      &0           \cr
0      &pI_2          \cr
\end{pmatrix},
$$
where $I_2$ is the $2\times 2$ identity matrix. \end{example}
\smallskip
The level function $u_r(n)$ that counts the
number of unitary factors of abelian groups of order $n$
whose conjugates have rank atmost $r$ is given by
$$ u_r(n)=\sum_{n=n_1^r n_2^{r-1}\dots n_r}
2^{\omega(n_1)+\omega(n_2)+\dots+\omega(n_r)},
$$
where $\omega(n)$ gives the number of distinct prime factors of $n$. We obtain
\begin{thm}
The associated zeta function \begin{gather*}
   \mathcal U_r(s)=\sum_{n=1}^\infty u_r(n) n^{-s} \\ \intertext{can be written as}  
\mathcal U_r(s)=\frac {\zeta^2(s)\zeta^2(2s)\zeta^2(3s)\cdots\zeta^2(rs)}
 {\zeta(2s)\zeta(4s)\zeta(6s)\cdots\zeta(2rs)}.
\end{gather*}
\end{thm}
This implies the following result on average orders:
\begin{cor} The number of unitary factors of finite abelian  groups of order
atmost $x$ whose conjugates have rank atmost $r$ is asymptotic to $K_rx\log
x$, where $K_r$ is a constant that depends on $r$. \end{cor}

In the infinite dimensional case, the study of the average order first occurs 
in
literature at the same time as direct factors and continues to be refined.
(See, eg. Calder{\'o}n \cite{Calderon} or Zhai \cite{Zhai}).

\section{Subgroup Classes}
It is easy to see that the number of classes of subgroups of $G$
corresponds to the number of divisors of $M(G)$ which are themselves in SNF.
Thus we are interested in the function 
$$ \sum_{M_1\mid M(G),\ M_1\ SNF}1=\sum_{H<G}1,$$
where for all isomorphic subgroups $H$, we count only once. For our purpose it
is enough to consider a $p$-group

$$ G_p \cong {\mathbb Z}/p^{a_1}{\mathbb Z} \oplus \cdots \oplus 
{\mathbb Z}/p^{a_1+\cdots+a_r}{\mathbb Z}$$ 
whose type is the partition $(a_1,a_1+a_2,\dots,a_1+\cdots+a_r),$  $a_i \ge 0$. We wish to
count all subgroups $H$ whose 
type is "less than" 
 the type of $G$, i.e. where the type of
$H$ is $(b_1,b_1+b_2,\dots,b_1+\cdots+b_r)$ with $\sum_{i=1}^j b_i \le \sum_{i=1}^j a_i$, for $1 \leq j\leq r$.

Thus we study the function 
\begin{gather} f_r(a_1,a_2,\dots,a_r)=\sum_{b_1=0}^{a_1}\sum_{b_2=0}^{a_1+a_2-b_1}\dots
\sum_{b_r=0}^{a_1+a_2+\dots+ a_r-b_1-b_2-\dots - b_{r-1}}1. \label{1} 
\end{gather}
 Collecting the exponents $r a_1+ \cdots+ a_r=k,$ we write 
\begin{gather}
 F_r(x) = \sum_{k=0}^\infty \alpha_r(k) x^k, \label{3} \\ \intertext{
 where}
 \alpha_r(k) = \sum_{r a_1+ \cdots  +a_r=k}f_r(a_1,a_2,\dots,a_r)  .  \label{4}
\end{gather}

Now let  $c_j= a_1+ \dots +a_j$ and $d_j=b_1+ \dots + b_j$ and let $\lambda _k=(c_1,c_2,\cdots ,c_r)$
and $\mu_m=(d_1,d_2,\cdots ,d_r)$  be   partitions  of $k$ and $m\le k$
respectively.
Then $\alpha _r(k)$ denotes the number of pairs of partitions
$\{(\lambda_k,\mu_m) \}$, such that $d_i \leq c_i$.

 We first obtain a recursion formula for the $\alpha_r(k)$ . To do so directly
 seems difficult, and it is convenient to introduce an additional parameter $j$.
 
\begin{defn} \label{arnj} Let $\alpha_r(k,j)$ be defined  as  $\alpha_r(k)$, but where the condition $0 \leq d_1 \leq \cdots \leq d_r$ is replaced by $-j \leq d_0 \leq d_1 \leq \cdots \leq d_r$, i.e. 
\begin{gather}
\alpha_r(k,j)=\sum_{\substack{c_1+ \cdots  +c_r=k \\ 0 \leq c_1 \leq \cdots \leq c_r}} \sum_{d_0=-j}^0\sum_{d_1=d_0}^{c_1}\sum_{d_2=d_1}^{c_2}\dots
\sum_{d_r=d_{r-1}}^{c_r}1.\label{8}
\end{gather}
\end{defn}

 \begin{prop} \label{anj} One has the recursion formula
 \begin{align*}
     &(i) &  \qquad  \alpha_0(k,j)&=\begin{cases}  j+1,  & k= 0, \\ 0, & k \geq 1. \end{cases} \\
     &(ii) &    \alpha_r(k,-j)&=0, \qquad j\ge 1.  \\
     &(iii) &      \alpha_r(k,j)&=\alpha_r(k,j-1) + \sum_{m=0}^{[k/r]} \alpha_{r-1}(k-mr,m+j),
       \qquad  (r \geq 1,j \geq 0).\\
     &(iv)& \alpha_r(k,0)&=\alpha_r(k). 
 \end{align*}
 \end{prop}
 \begin{proof}
  Equations $(i), (ii)$ and  $(iv)$ follow  from  
Definition \ref{arnj}.  

For  $(iii)$,
   the  condition $-j \leq d_0$ in Definition 4 can be divided into two cases, 
either $-(j-1) \leq d_0$ from which we get the contribution
$\alpha_r(k,j-1),$ in  $(iii)$ ; 
or $d_0=-j$ in which case we let $c_1=m$ and denote
$$
     c_i'=c_{i+1}-m, \qquad \text{and} \qquad d_i'=d_{i+1}-m. 
$$ 
It is clear that
$ (0,c_1', c_2', \cdots \leq c_{r-1}')$ is a partition of $k-mr$., 
while  $d_0' \leq 0, \qquad -m-j \leq d_0' \leq d_1' \leq \dots \leq d_{r-1}',
\qquad  \text{ and } \qquad d_j' \leq c_j'.$ 
 We thus get a contribution 
\begin{gather*}  \alpha_{r-1}(k-mr,m+j) 
\end{gather*}
 for each $0 \leq m \leq [k/r]$, and this concludes the proof of $(iii)$.
 \end{proof}
We introduce the notation 
$\alpha(k)= \alpha_{k}(k).$

Since the number of partitions of $k$ into $k$ non-zero parts is equal to the
number of partitions of $k$ into $r$ parts, $r>k$, we get,
\begin{gather} \label{alphadef}
\alpha_r(k)= \alpha (k).\qquad (r\ge k) \\ \intertext{Thus,} \notag
 \lim_{r \to \infty} \alpha_r(k)=\alpha(k)
\end{gather}
We now obtain some simple bounds in terms of the  classical partition function.
\begin{lem} \label{ppp}
  We have that 
  \begin{gather*}
    q(k) \leq \alpha(k) \leq kq(k)^2, 
   \end{gather*} 
    where $q(n)$ denotes the total number of partitions of the $n$. 
\end{lem}

\begin{proof}
The lower bound follows from the definition. Now  
let $\lambda _k=(c_1,c_2,\cdots ,c_k)$
and $\lambda _m=(d_1,d_2,\cdots ,d_k)$  be   partitions  of $k$ and $m\le k$
respectively. If $d_1\ne 0, \lambda _m $ can be completed to
$(1,\cdots,1,d_1,d_2,\cdots ,d_k)$, a partition of $k$. There are atmost
$q(k)$ of such $\lambda _m $.

If $d_1=0,d_2\ne 0$, we can complete $\lambda _m $ to get a partition of
$k-c_1$. In this way, for each $\lambda _m $ we have atmost
$$ q(k)+q(k-c_1)+\cdots+q(k-c_k)\le kq(k)$$
possibilities. Since there are $q(k)$ ways to write $\lambda _k $, the upper
bound follows.
\end{proof}

We will now turn our attention to the generating functions
\begin{defn} \label{Fdef} Let $\alpha_r(k), \alpha_r(k,j),$ and $\alpha(k)$ be
  as above. We define the generating functions to be  
\begin{align*}
 F_r(x)&=  \sum_{k=0}^\infty \alpha_r(k) x^k, \\
 F_r(x,y)&= \sum_{k,j=0}^\infty \alpha_r(k,j) x^k y^j, \\ \intertext{and}
  F(x)&=  \sum_{k=0}^\infty \alpha(k) x^k. 
\end{align*}
\end{defn}

These definitions are consistent, since

 \begin{gather}\label{frx}
   F_r(x)=F_r(x,0),
 \end{gather}
the only term that survives with $y=0$ is when $j=0$.

\begin{lem} \label{u13}
 With $F$ and $F_r$ defined as above one has that $F(x)$ and $F_r(x)$ are analytic
 functions in the unit disc with integer power series coefficients such that $F(0)=F_r(0)=1$. Furthermore the function $F$ 
satisfies the inequality
 \begin{gather*}
    F(x) \geq \frac 1 {\prod_{k=1}^\infty (1-x^k)}. \qquad \qquad (0<x<1) 
 \end{gather*}
\end{lem}

\begin{proof}
The power series coefficients of $F_r$ and $F$ are integers since they are counting functions and 
by Proposition 2 $(i)$ and $(iv)$ and eq. \eqref{alphadef}, we have that $\alpha_r(0)=1$ and $\alpha(0)=1$, which implies $F_r(0)=F(0)=1$. 
By the well known generating function for the classical partition function  
 \begin{gather} \label{o112}
  \sum_{k=0}^\infty q(k) x^k= \frac 1 {\prod_{k=1}^\infty (1-x^k)}, \qquad \qquad (0<x<1)
  \end{gather}
 and the lower bound in Lemma \ref{ppp}
 \begin{gather*}
   q(k) \leq \alpha(k),
 \end{gather*}
this gives us the lower bound in Lemma \ref{u13}. Equation \eqref{o112} also
implies that the 
generating function of the partition function is analytic in the unit disc,
and hence the 
classical partition function $q(n)$ is of subexponential order. This implies
that   $k (q(k))^2$ is 
of subexponential order and by the upper bound in Lemma \ref{u13}, so is
$\alpha(n)$, and also $\alpha_r(n)$ since $0 \leq \alpha_r(n)\leq \alpha(n)$.
This  proves that $F$ and $F_r$ are analytic in the unit disc. 
\end{proof}

We will now see how the recursion formula for the coefficients
$\alpha_r(k,j)$ in Proposition \ref{anj} 
 yield a recursion formula for its
generating function $F_r(x,y)$.

\begin{prop} \label{frec}  One has the following recursion formula for the function $F_r(x,y)$:
\begin{align*}
   &(i) \qquad  &  F_0(x,y)&=\frac 1 {(1-y)^2}, \\   
   &(ii) &  F_r(x,y)&=\frac  {F_{r-1}(x,x^r)x^r  -  F_{r-1}(x,y) y}{(x^r-y)(1-y)}.
\end{align*}
\end{prop}
\begin{proof} From Proposition \ref{anj} $(iii)$ it is clear that
$$
  F_r(x,y)(1-y) = \sum_{k,j=0}^\infty \sum_{m=0}^{[k/r]} \alpha_{r-1}(k-mr,m+j) x^n y^j,
$$
and by replacing $k$  and $j$ suitably, we get
\begin{align*}
 F_r(x,y)- F_r(x,y) y 
  &=  \sum_{k,j=0}^\infty \sum_{m=0}^j x^{k+mr} y^{j-m} \alpha_{r-1}(k,j),   \\ 
 &= \sum_{k,j=0}^\infty x^{k} \left(\frac{x^{r(j+1)}- y^{j+1}}{x^r-y} \right) \alpha_{r-1}(k,j),  \\ 
    &= \frac  {F_{r-1}(x,x^r)x^r -   F_{r-1}(x,y) y}{x^r- y}.
\end{align*}
Proposition  \ref{anj} $(i)$ and Definition 5 give 
\begin{align*}
  F_0(x,y)&=\sum_{j=0}^\infty (j+1) y^j, \\&= \frac 1 {(1-y)^2}.
\end{align*}

\end{proof}
We first notice that
\begin{lem} \label{vw1} The function $F$ satisfies the identity.
\begin{gather*}
 F_r(x)=  F_{r-1}(x,x^r)
\end{gather*}
\end{lem}
\begin{proof}
 This follows from eq. \eqref{frx} 
 and  Proposition  \ref{frec} $(ii)$. 
\end{proof}

We now calculate the first few values of $F_r(x)$.

\begin{lem} \label{Fex}
  One has that
  \begin{align*}
    &(i) \qquad & F_0(x)&=1, \\
    &(ii) & F_1(x)&=(1-x)^{-2}, \\ \intertext{and}
    &(iii) & F_2(x)&=(1-x)^{-2}(1-x^2)^{-3} (1-x^3).
  \end{align*}
\end{lem}
\begin{proof}
 Putting $y=0$ in Proposition \ref{frec} $(i)$  
    proves $(i)$. From Lemma \ref{vw1} we get that
  \begin{gather*}
    F_1(x)=F_0(x,x)=\frac{1}{(1-x)^2},
  \end{gather*}
  which proves $(ii)$.
 By Proposition \ref{frec} $(ii)$ we get that
   \begin{gather} \label{ww} \begin{split}
     F_1(x,y)&=\frac  {F_{0}(x,x)x  -  F_{0}(x,y) y}{(x-y)(1-y)} \\ 
             &=\frac{xy-1}{(x-1)^2 (y-1)^3}. \end{split}
   \end{gather}
    By Lemma \ref{vw1} we get
   \begin{gather*}
      F_2(x)=F_{1}(x,x^2)=      \frac{x^3-1}{(x-1)^2(x^2-1)^3},
    \end{gather*}
    which proves $(iii)$.
\end{proof}

We further prove:
\begin{prop} The  functions $F_r(x,y)$ and $F_r(x)$ satisfy the
functional equations  \label{cor3}
\begin{align*} 
 &(i) \qquad &   F_r(x,y)&=x^{-r(r+1)/2}y^{-2}F_r \left(\frac 1 x, \frac 1 y \right), 
 \\
  &(ii) &  F_r(x)&=x^{-r(r+3)/2} F_r \left(\frac 1 x \right).
\end{align*}
\end{prop}
\begin{proof} We will use induction to  prove $(i)$.  Let $r=0$.
Proposition \ref{frec} $(i)$ gives us that 
\begin{align*}
   F_0(x,y) &= (1-y)^{-2},\\  &=y^{-2}\left(1-\frac 1 y \right)^{-2}, \\
   &= y^{-2} F_0 \left(\frac 1 x, \frac 1 y \right), \\ &=
 x^{-0\cdot(0+1)/2} y^{-2} F_0 \left(\frac 1 x,\frac 1 y \right),
\end{align*}
 and hence $(i)$ is true for $r=0$. We now assume that $(i)$ is true for
 $r=k$. 
By using  Proposition \eqref{frec} $(ii)$ we now obtain
\begin{align*}
    F_{k+1}(x,y)&=\frac{F_{k}(x,x^{k+1})x^{k+1}  -  F_{k}(x,y)
      y}{(x^{k+1}-y)(1-y)},
 \\ \intertext{Using the functional equation $(i)$ for $r=k$ it equals} 
 &=\frac{x^{-k(k+1)/2} x^{-2(k+1)} F_{k} (x^{-1},x^{-k-1})x^{k+1}  - 
   x^{-k(k+1)/2} y^{-2} F_{k}(x^{-1},y^{-1}) y}{(x^{k+1}-y)(1-y)}, \\ 
 &=x^{-(k+1)(k+2)/2} y^{-2} \frac{F_{k} (x^{-1},x^{-k-1})x^{-k-1}  - 
    F_{k}(x^{-1},y^{-1}) y^{-1}}{(x^{-k-1}-y^{-1})(1-y^{-1})}, \\  &=
 x^{-(k+1)(k+2)/2} y^{-2} F_{k+1}(x^{-1},y^{-1}),
  \end{align*}
which finishes the proof for $(i)$.
We will now prove $(ii)$. By Lemma \ref{vw1} we have that
\begin{gather*}
 F_r(x) = F_{r-1}(x,x^r), \\ \intertext{The functional equation with $y=x^r$ gives}
   x^{-(r-1) r/2} x^{-2r} 
   F_{r-1} \left( \frac 1 x, \frac 1 {x^r} \right),
 \\ 
  = x^{-r(r+3)/2} F_{r}\left( \frac 1 x \right).
\end{gather*}
\end{proof}

\section{The zeta function of subgroup classes}

\subsection{The zeta function of subgroup classes of bounded rank}
We will now introduce the zeta function.  We use the notation $A_r(n)$ for the multiplicative 
function that counts the number of isomorphic subgroup classes
of abelian groups of rank less or equal to $r$ and order $n$. If $n$ is a prime power then
$$A_r(p^k)=\alpha_r(k).$$
The corresponding zeta function is given by
\begin{defn} \label{zetar}Let \begin{gather*}
    Z_r(s) = \sum_{n=1}^\infty A_r(n) n^{-s}=\prod_{p \text{ prime}} F_r(p^{-s}).
   \end{gather*}\end{defn}
We then get

\begin{lem} \label{lem99}
 One has that
 \begin{gather*}
   Z_r(s)= \prod_{m=1}^k \zeta(ms)^{\beta_{r}(m)} G_{k,r}(s),
  \end{gather*} where   $G_{k,r}(s)$
  is a Dirichlet series without real zeroes and absolutely convergent for $\Re(s) >1/(k+1)$.
  Furthermore the $\beta_{r}(m)$'s are integers.
\end{lem}

\begin{proof}
By Lemma 2 the result follows from Dahlquist \cite{D} with
 \begin{gather} \label{betardef}
 \beta_{r}(m)=  \sum_{d|m} \mu \left(\frac m d \right) \frac d m B_r(d),
\\ \intertext{where $B_r$ comes from the expression}
 \notag \log F_r(x)  = \sum_{m=1}^{\infty}B_r(m)x^{m}. \\ \intertext{Thus} \notag
 \begin{split}
    G_{k,r}(s) &= \prod_{p \text{ prime}} F_r(p^{-s}) \prod_{m=1}^k (1-p^{-ms})^{-\beta_{r}(m)} \\ &=  
\prod_{m=1}^k \zeta(ms)^{\beta_{r}(m)} Z_r(s)
\end{split}
\end{gather}
is without real zeroes and absolutely convergent for $\Re(s)>1/(k+1)$. 
\end{proof}
 With the notation $ \beta(m)=\beta_{m}(m)$ we get from  eq. \eqref{alphadef} and the definition of the $\beta_r(m)$ eq. \eqref{betardef} that
\begin{gather}
    \beta(m)=\beta_{r}(m). \qquad (r \geq m) \label{betadef}
\end{gather}
By using Definition 6, Lemma \ref{lem99}, eq \eqref{betadef} and Lemma 4 we obtain 
\begin{thm}
   One has that
   \begin{align*}
     Z_1(s)&=\zeta^2(s), \\ Z_2(s)&=\frac{\zeta^2(s) 
\zeta^3(2s)}{\zeta(3s)}, \\ \intertext{and}  \quad Z_k(s)&= \zeta^2(s) \zeta^3(2s) \zeta^3(3s) G_{3,k}(s), \qquad (k \geq 3)
    \end{align*}
    where $G_{3,k}(s)$ is a Dirichlet series without real zeroes and absolutely convergent for $\Re(s)>1/4$.
 \end{thm}
This implies the following result on average orders:
 \begin{cor}  The number of subgroup classes of abelian groups of rank less than or equal to $r$ and order atmost $x$ is asymptotic to 
$C_r x \log x$ where $C_r$ is a constant depending on $r$.
\end{cor}

\subsection{The zeta function of subgroup classes of arbitrary rank}
 
  We now let $A(n)=\lim_{r \to \infty} A_r(n)$, $Z(s)=\lim_{r \to \infty} Z_r(s)$ and proceed as in the last
  section. We obtain with Lemma 2 and Dahlquist's theorem that
 \begin{lem} \label{lem99c}
 One has that
 \begin{gather*}
   Z(s)= \prod_{m=1}^k \zeta(ms)^{\beta(m)} G_{k}(s), \end{gather*} where
  $G_{k}(s)$ 
  is a Dirichlet series without real zeroes and absolutely convergent for $\Re(s) >1/(k+1)$, and the $\beta(m)$'s are integers defined by eq. \eqref{betadef}.
\end{lem}
which implies
   \begin{thm}
 Let $A(n)$ denote the number of subgroup classes of abelian groups of order $n$, and let 
  \begin{gather*}
    Z(s) = \sum_{n=1}^\infty A(n) n^{-s}
   \end{gather*}
   denote the corresponding zeta function. Then 
   \begin{gather*}
     Z(s)= \zeta^2(s) \zeta^3(2s) \zeta^3(3s) G_{3}(s) 
   \end{gather*}
   where $G_{3}(s)$ is a Dirichlet series without real zeroes and absolutely convergent for $\Re(s)>1/4$. Furthermore $Z(s)$ is a meromorphic function for $\Re(s)>0$ and the 
imaginary axis is the natural boundary for  $Z(s)$.
\end{thm}
\begin{proof} That $\Re(s)=0$ is the natural
 boundary follows from Lemma 2 which shows that 
  $F$ can not be written as a finite product
$\prod_{j=1}^k (1-x^j)^{m_j}. \qquad (m_j \in \mathbb Z)$ 
  and  from  Dahlquist's theorem  the function cannot be meromorphically
  continued beyond the line $\Re(s)=0$. \end{proof}

As before we get an estimate:
\begin{cor} The number of subgroup classes of abelian groups of order atmost $x$ is asymptotic to 
$C x \log x$ where $C$ is a constant.
\end{cor}

\begin{rem}
  The constant $C$ in Corollary 4 can be calculated by 
  \begin{gather*}
    C=\prod_{m=2}^\infty \zeta(m)^{\beta(m)}  =13.1854452968422695.
  \end{gather*}
\end{rem}

We notice that $\beta(m) \geq 0$ for $m=1,\dots,12$ and for $k=13$ Lemma \ref{lem99c} allows us to see (with the Mathematica program in 
the Appendix)  
\begin{rem}
 One has that 
\begin{multline*}
      Z(s)= \zeta^2(s) \zeta^3(2s) \zeta^3(3s) \zeta^4(4s) \zeta^4(5s) 
      \zeta^6(6s) \times \\ \times \zeta(7s) \zeta^4(8s)
 \zeta^6(9s) \zeta^2(10s)  \zeta^{12}(12s) (\zeta(13s))^{-1} G_{13}(s),
   \end{multline*}
  where $G_{13}(s)$ is a Dirichlet
 series absolutely convergent for $\Re(s)>1/14$
  In particular this means that there are no poles for $Z(s)$ for $\Re(s)>1/13$
 except for $s=1/k$, $k$ integer. Under the Riemann hypothesis, we can use Lemma 6 with $k=25$ to see that there are no poles for $Z(s)$ for $\Re(s)>1/26$ except for $s=1/k$, $k$ integer.
\end{rem}
For more calculations of $\beta(m)$ see the appendix where we use a Mathematica program to calculate its value for the first $100$ values of $m$. 
\subsection{A question on the zeta function of subgroup classes of bounded rank}
 It seems likely 
 that $r=1$ and $2$ are the only cases where $Z_r(s)$ has a meromorphic
 continuation to the entire complex plane. 
For higher rank we expect that

{\it $Z_r(s)$ is a meromorphic function for $\Re(s)>0$ with the 
imaginary axis as its 

natural boundary.}

\noindent That $Z_r(s)$ has a meromorphic continuation to $\Re(s)>0$ of course follows from Lemma 5. To prove
that $\Re(s)=0$ is the natural boundary we would have to show that $F_r(x)$ cannot be written as a finite product
\begin{gather*}
  F_r(x)= \prod_{j=1}^k (1-x^j)^{m_j}. \qquad (m_j \in \mathbb Z)
\end{gather*}
Thus it is sufficient to show that not all roots of the numerator of 
$F_r(x)$  lie on the unit circle which is easily verifiable 
by numerical calculation for small values of $r$. 
As an example, let 
$r=3$. We can use the explicit form of $F_3(x)$ as given in the appendix, 
and  see  that the equation
$1+2x+2x^2+x^3+x^4=0$ has a root 
$x= -0.621744 + 0.440597 i$ with an absolute value $0.762031$.

Such calculations are not obvious for a general $r$.

\section{Appendix}

\subsection{Tables}

\begin{center}
\begin{tabular}{|l|cccccccccccccccc|} \hline  
\multicolumn{17}{|c|}{ Table of $\beta_{r}(m)$}\\\hline
\hline  
  $r \backslash m$ & 1 & 2 & 3 & 4 & 5 & 6 & 7 & 8 & 9 & 10 & 11 & 12 & 13 & 14 & 15 &16  \\ \hline  
  1 & 2 & 0 & 0 & 0 & 0 & 0 & 0 & 0 & 0 & 0 & 0 & 0 & 0 & 0 & 0 & 0 \\
  2 &  2 & 3 & -1 & 0 & 0 & 0 & 0 & 0 & 0 & 0 & 0 & 0 & 0 & 0 & 0 & 0 \\
  3 &  2 & 3 & 3 & -1 & -3 & 3 & 0 & -2 & 0 & 0 & 6 & -5 & -6 & 6 & 1 & 9 \\
  4 &   2 & 3 & 3 & 4 & -4 & 0 & -3 & 2 & 9 & -10 & -2 & -2 & 3 & 25 & -25 & 11 \\
  5 &  2 & 3 & 3 & 4 & 2 & -1 & -6 & -1 & 3 & 4 & 2 & 3 & -12 & -22 & 15 & 26 \\
  6 &  2 & 3 & 3 & 4 & 2 & 6 & -7 & -4 & 0 & -2 & 1 & 13 & 10 & -21 & -4 & -9 \\
  7 &  2 & 3 & 3 & 4 & 2 & 6 & 1 & -5 & -3 & -5 & -5 & 12 & -1 & 8 & 17 & -9 \\
  8 &  2 & 3 & 3 & 4 & 2 & 6 & 1 & 4 & -4 & -8 & -8 & 6 & -2 & -3 & 18 & 20 \\
  9 &  2 & 3 & 3 & 4 & 2 & 6 & 1 & 4 & 6 & -9 & -11 & 3 & -8 & -4 & 7 & 21 \\
  10 &  2 & 3 & 3 & 4 & 2 & 6 & 1 & 4 & 6 & 2 & -12 & 0 & -11 & -10 & 6 & 10 \\
  11 &  2 & 3 & 3 & 4 & 2 & 6 & 1 & 4 & 6 & 2 & 0 & -1 & -14 & -13 & 0 & 9 \\
  12 &  2 & 3 & 3 & 4 & 2 & 6 & 1 & 4 & 6 & 2 & 0 & 12 & -15 & -16 & -3 & 3 \\
  13 &  2 & 3 & 3 & 4 & 2 & 6 & 1 & 4 & 6 & 2 & 0 & 12 & -1 & -17 & -6 & 0 \\
  14 &  2 & 3 & 3 & 4 & 2 & 6 & 1 & 4 & 6 & 2 & 0 & 12 & -1 & -2 & -7 & -3 \\
  15 &  2 & 3 & 3 & 4 & 2 & 6 & 1 & 4 & 6 & 2 & 0 & 12 & -1 & -2 & 9 & -4 \\
  16 &  2 & 3 & 3 & 4 & 2 & 6 & 1 & 4 & 6 & 2 & 0 & 12 & -1 & -2 & 9 & 13 \\ \hline
\end{tabular}
\end{center}

For the case of $\beta(m)$ we have the first 100 values: $\beta(1),\ldots,\beta(100)=$ $2,$ $3,$ $3,$ $4,$ $2,$ $6,$ $1,$ $4,$ $6,$ $2,$ $0,$ $12,$ $-1,$ $-2,$ $9,$ $13,$ $0,$ $-16,$ $6,$ $44,$
 $-25,$ $-6,$ $16,$ $-19,$ $52,$ $-21,$ $52,$ $-103,$ $-140,$ $505,$ $-203,$ $-381,$ $286,$ $88,$ $185,$ $-751,$ $564,$ $1015,$ $-2304,$ $1007,$ $1876,$ $-3432,$ $1177,$ $3665,$ $-1582,$ $-6119,$ $-2558,$ $21792,$ $-7745,$ $-34936,$ $40625,$ $4248,$ $-34948,$ $19176,$ $-361,$ $30668,$ $-104511,$ $81530,$ $196372,$ $-457425,$ $148194,$ $497952,$ $-459549,$ $-261973,$ $163195,$ $1015378,$ $-808365,$ $-1895457,$ $3467251,$ $-483924,$ $-4303223,$ $4879125,$ $535751,$ $-4991204,$ $-634706,$ $7525125,$ $7551313,$ $-31891601,$ $10748210,$ $53812155,$ $-65485063,$ $-16092470,$ $74326725,$ $-12532241,$ $-51617936,$ $-53724288,$ $219808510,$ $-91391802,$ $-403364243,$ $663219049,$ $-34111622,$ \\ $-908127917,$ $645523852,$ $643374980,$ $-402710571,$ $-1752056675,$ $1385661209,$ \\ $3632651805,$  $-6352203822,$ $-317398867$

\subsection{Mathematica Programs}
We can now implement the recursion formulae in Mathematica. 

\begin{verbatim}
Clear[F];
F[0,x_,y_] := F[0,x,y]=1/(1-y)^2;
F[r_,x_,y_]:= F[r,x,y] = 
   Factor[(x^r F[r-1,x,x^r]-F[r-1,x,y] y)/((x^r-y)(1-y))];
\end{verbatim}
Now using
\begin{verbatim}
F_r(x)=F[r,x,0]
\end{verbatim} 
 gives  us for $r=0,1,2,3$ that
\begin{align*} F_0(x)&=1, \\
   F_1(x)&= \frac 1 {{\left( -1 + x \right) }^{2}},
 \\
  F_2(x)&= \frac{1 + x + x^2}{{\left( -1 + x \right) }^4 {\left( 1 + x \right) }^3},
 \\ \intertext{and}
  F_3(x)&=\frac{\left( 1 + 2 x + 2 x^2 + x^3 + x^4 \right)  
     \left( 1 + x + 2 x^2 + 2 x^3 + x^4 \right) }{{\left( -1 + x \right) }^
      6 {\left( 1 + x \right) }^3 {\left( 1 + x + x^2 \right) }^4}.
 \end{align*}
The program above is  slow if we just want to calculate the lower order coefficients  $\alpha_r(n)$ in the power series of $F_r(x)$, so
the program we used to calculate the Table of $\beta_{r}(m)$  is the following (in this case with $m=16$) :
\begin{verbatim}
Clear[a];
a[r_,n_,-1]:=a[r,n,-1]= 0;
a[0,n_,j_]:= a[0,n,j] = If[n==0,j+1,0];
a[r_,n_,j_]:=a[r,n,j]=a[r,n,j-1]+Sum[a[r-1,n-m r,m+j],
             {m,0,Floor[n/r]}];
a[r_,n_]:=a[r,n]=a[r,n,0];

m=16; Table[aa=Series[Log[Sum[a[r,j]*x^j,{j,0,m}]],{x,0,m}][[3]];   
Table[ b=Divisors[n];
Sum[MoebiusMu[b[[j]]]/b[[j]] aa[[n/b[[j]]]],{j,1,Length[b]}],
{n,1,m}],{r,1,m}]
\end{verbatim}
A similar program is used to calculate the $\beta(m)$ for $m=1,\dots 100$.

\bibliographystyle{alpha}

\end{document}